\documentclass[a4paper]{article}
\usepackage{graphicx}
\usepackage{amssymb}
\usepackage{amsmath}
\usepackage{amsfonts}
\usepackage{varioref}
\usepackage[bookmarksnumbered]{hyperref}
\usepackage[numbers]{natbib}

\title{ Domain decomposition methods for problems of unilateral contact between elastic bodies with nonlinear Winkler covers
        \footnote{This work was partially supported by Grant 23-08-12 of National Academy of Sciences of Ukraine}}

\author{Ihor~I.~Prokopyshyn
\footnote{Pidstryhach Institute for Applied Problems of Mechanics
and Mathematics, Naukova 3-b, Lviv, 79060, Ukraine,
ihor84@gmail.com, Corresponding author }\and Ivan~I.~Dyyak
\footnote{Ivan Franko National University of Lviv, Universytetska
1, Lviv, 79000, Ukraine, dyyak@franko.lviv.ua} \and
 Rostyslav~M.~Martynyak \footnote{Pidstryhach Institute for Applied Problems of Mechanics and Mathematics,
Naukova 3-b, Lviv, 79060, Ukraine, labmtd@iapmm.lviv.ua } \and
Ivan~A.~Prokopyshyn \footnote {Ivan Franko National University of
Lviv, Universytetska 1, Lviv, 79000, Ukraine, lviv.pi@gmail.com}}

\begin{document}
\maketitle

\begin{abstract}
%% Text of abstract
In this paper we propose on continuous level a class of domain decomposition methods of
Robin--Robin type to solve the problems of unilateral contact between elastic bodies
with nonlinear Winkler covers. These methods are based on abstract nonstationary
iterative algorithms for nonlinear variational equations in reflexive Banach spaces. We
also provide numerical investigations of obtained methods using finite element
approximations.


\end{abstract}

\section{Introduction} \label{prokopyshyni_contrib:1}
Thin covers from another material are often applied in engineering
to improve the functional properties of the surfaces of components
of machines and structures. On the other hand, thin covers with
certain mechanical properties are used for modeling of real
microstructure of the surfaces, adhesion and glue bondings
\cite{prokopyshyni_contrib_Goryacheva1998,
prokopyshyni_contrib_Suquet1988,
prokopyshyni_contrib_ContMech2001}.

The classical methods for solution of contact problems for bodies
with thin covers are grounded on integral equations and are
reviewed in work \cite{prokopyshyni_contrib_ContMech2001}.
Nowadays, one of the most effective numerical methods for such
contact problems are methods, based on variational formulations
and finite element approximations.

Efficient approach for solution of multibody contact problems is the use of domain
decomposition methods (DDMs). Many DDMs for contact problems without covers are
obtained on discrete level
\cite{prokopyshyni_contrib_Dostal2010,prokopyshyni_contrib_Wohlmuth2011}. Among DDMs,
proposed on continuous level for contact problems without covers are methods presented
in
\cite{prokopyshyni_contrib_Bayada2004,prokopyshyni_contrib_Koko2009,prokopyshyni_contrib_Sassi2008}.
Domain decomposition methods for solution of problem of ideal contact between two
bodies, connected through nonlinear Winkler layer are proposed in
\cite{prokopyshyni_contrib_Bresch2004, prokopyshyni_contrib_Koko2008}. These methods
are based on saddle-point formulation and conjugate gradient methods.

In current contribution we consider the problem of unilateral
contact between bodies with nonlinear Winkler covers. We give
variational formulations of this problem in the form of nonlinear
variational inequality on convex set and variational equation in
the whole space, and present theorems about existence and
uniqueness of their solution. Furthermore, we propose on
continuous level a class of parallel domain decomposition methods
for solving the nonlinear variational equation, which corresponds
to original contact problem. In each iteration of these methods we
have to solve in a parallel way linear variational equations in
separate bodies, which are equivalent in a weak sense to linear
elasticity problems with Robin boundary conditions on possible
contact areas. These DDMs are based on abstract nonstationary
iterative methods for variational equations in Banach spaces. They
are the generalization of domain decomposition methods, proposed
by us earlier in \cite{prokopyshyni_contrib_Dyyak2010b,
prokopyshyni_contrib_Dyyak2012,prokopyshyni_contrib_Prokopyshyn2008}
for unilateral contact problems without covers. Some particular
cases of proposed DDMs can be viewed as a modification of
semismooth Newton method
\cite{prokopyshyni_contrib_Hintermuller2003}. The numerical
analysis of obtained DDMs is made for plane contact problems using
finite element approximations.

\section{Statement of the problem} \label{prokopyshyni_contrib:2}
% Always give a unique label
% and use \ref{<label>} for cross-references
% and \cite{<label>} for bibliographic references
% use \sectionmark{}
% to alter or adjust the section heading in the running head

Consider a unilateral contact of $N$ elastic bodies $\Omega _{\alpha } \subset {\mathbb
R}^{3} $ with piecewise smooth boundaries $\Gamma _{\alpha }$, $\alpha =1,2,...,N$
(Fig.~1a). Suppose that %on
across each contact surface there is a nonlinear Winkler layer.
Denote $\Omega =\bigcup_{\alpha =1}^{N}\Omega_{\alpha}$.

\begin{figure}[h]
%\sidecaption

\includegraphics[scale=1]{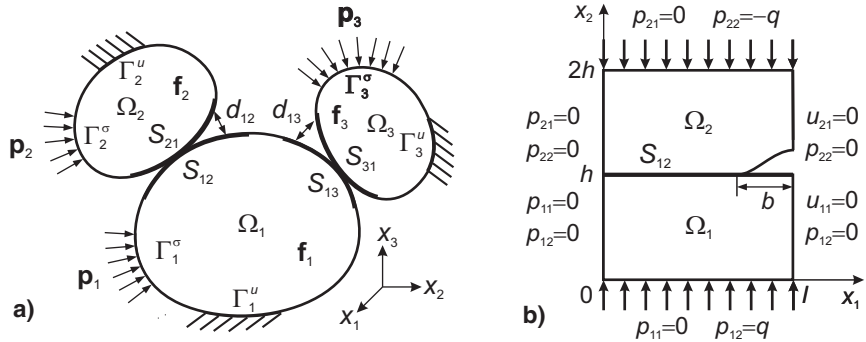}
\caption{Unilateral contact between several elastic bodies through
nonlinear Winkler layers}
\label{author_program_element_fig:1}       % Give a unique label
\end{figure}

A stress-strain state in point ${\bf x}=(x_{1}
,x_{2},x_{3})^{\top}$ of each solid $\Omega _{\alpha }$ is
described by the displacement vector ${\bf u}_{\,\alpha}=u_{\alpha
\, i} \, {\bf e}_{i}\,$, the tensor of strains $
{\hat{\pmb{\varepsilon}}}_{\alpha }=\varepsilon _{\alpha \,
ij}\,{\bf e}_{i} \, {\bf e}_{j} $ and the tensor of stresses
$\hat{{\pmb \sigma }}_{\alpha } =\sigma _{\alpha \, ij} \, {\bf
e}_{i} \, {\bf e}_{j}$. These quantities satisfy the following
relations:
\begin{equation} \label{prokopyshyni_contrib__1_}
\sum _{j=1}^{3} \frac{\partial \sigma _{\alpha \, ij}({\bf
x})}{\partial x_{j}} \, +f_{\alpha \, i}({\bf x}) =0\,, \,\,\,\,
{\bf x} \in {\Omega}_{\alpha}\,, \,\,\,\, i=1,2,3\,,
\end{equation}
\begin{equation} \label{prokopyshyni_contrib__1b_}
\sigma _{\alpha \,  ij}({\bf x}) = \sum _{k,l=1}^{3} C_{\alpha \,
ijkl}({\bf x}) \, \varepsilon _{\alpha \, kl}({\bf x})\,, \,\,\,\,
{\bf x} \in {\Omega}_{\alpha}\,, \,\,\,\, i,j=1,2,3\,,
\end{equation}
\begin{equation} \label{prokopyshyni_contrib__2_}
\varepsilon _{\alpha \, ij}({\bf x}) =\frac{1}{2}
\left(\frac{\partial u_{\alpha \, i}({\bf x})}{\partial x_{j} }
+\frac{\partial u_{\alpha \, j}({\bf x})}{\partial x_{i} }
\right)\,,\,\,\,\, {\bf x} \in {\Omega}_{\alpha}, \,\,\,\,
i,j=1,2,3\,,
\end{equation}
where $f_{\alpha \, i}$ are the components of volume forces vector
${\,\bf f}_{\alpha } =f_{\alpha \, i} \,{\bf e}_{i} $, and
$C_{\alpha \, ijkl}$ are symmetric elasticity constants, which are
bounded in the following sense:
\begin{equation} \label{prokopyshyni_contrib__3_}
\left(\exists b_{\alpha},c_{\alpha}>0 \right)\left(\forall {\bf x}
\right)\left\{b_{\alpha}\sum _{i,j=1}^{3}\varepsilon _{\alpha
ij}^{2}
 \le \sum _{i,j,k,l=1}^{3}C_{\alpha ijkl}\varepsilon _{\alpha ij}
\varepsilon _{\alpha kl}  \le c_{\alpha}\sum
_{k,l=1}^{3}\varepsilon _{\alpha kl}^{2}\right\}.
\end{equation}

Introduce on boundary $\Gamma _{\alpha }$ a local orthonormal
coordinate system ${\pmb \xi }_{\alpha } ,\, { \pmb \eta }_{\alpha
} ,\, {\bf n}_{\alpha }$, where ${\bf n}_{\alpha }$ is an outer
unit normal. Then the vectors of displacements and stresses on
$\Gamma_{\alpha}$ can be written in the following way: ${\bf
u}_{\,\alpha } =u_{\alpha \, \xi } \, { \pmb{\xi} }_{\alpha }
+u_{\alpha \eta } \, {\pmb{\eta} }_{\alpha } +u_{\alpha n} \,
\bf{n}_{\,\alpha } , \, \,
 {\pmb{ \sigma} }_{\alpha } =\hat{{\pmb{ \sigma }}}_{\alpha }
\cdot n_{\,\alpha } =\sigma _{\alpha \xi } \, {\pmb \xi }_{\alpha }
 +\sigma _{\alpha \eta } \, {\pmb \eta }_{\alpha } +\sigma _{\alpha \it n}
 \, n_{\,\alpha }\,.$

Suppose, that the boundary $\Gamma _{\alpha } $ consists of three
disjoint parts:\\ $\Gamma _{\alpha } =\Gamma _{\alpha }^{u}
\bigcup \Gamma _{\alpha }^{\sigma } \bigcup S_{\alpha } $, $\Gamma
_{\alpha }^{u} =\overline{\Gamma _{\alpha }^{u} }$, $\Gamma
_{\alpha }^{u} \ne \emptyset $, $S_{\alpha } \ne \emptyset $. On
the part $\Gamma_{\alpha }^{u}$ homogenous Dirichlet boundary
conditions are prescribed, and on the part $\Gamma _{\alpha
}^{\sigma}$ we consider Neumann boundary conditions:
\begin{equation} \label{prokopyshyni_contrib__4_}
{\bf u}_{\,\alpha } ({\bf x})=0, \,\,\, {\bf x}\in \Gamma _{\alpha}^{u}
\,;\,\,\,\,\,{\pmb \sigma }_{\alpha } ({\bf x})={\bf p}_{\alpha } ({\bf x}), \,\,\,
{\bf x}\in \Gamma _{\alpha }^{\sigma } \, .
\end{equation}

The part $S_{\alpha } =\bigcup _{\beta \in B_{\alpha } }S_{\alpha
\beta }$, $\bigcap_{\beta \in B_{\alpha } }S_{\alpha \beta } =
\emptyset $ is the possible contact area of body $\Omega_{\alpha}$
with the other bodies. Here $S_{\alpha \beta } $ is the possible
unilateral contact area of body $\Omega _{\alpha } $ with body
$\Omega _{\beta } $, and $B_{\alpha } \subset
\left\{1,2,...,N\right\}$ is the set of the indices of all bodies
in contact with body $\Omega _{\alpha } $. We assume that the
surfaces $S_{\alpha \beta } \subset \Gamma _{\alpha } $ and
$S_{\beta \alpha } \subset \Gamma _{\beta } $ are sufficiently
close ($S_{\alpha \beta } \approx S_{\beta \alpha } $), and  ${\bf
n}_{\,\alpha } ({\bf x})\approx -{\bf n}_{\,\beta } ({\bf x}')$,
${\bf x}\in S_{\alpha \beta } $,  ${\bf x}'=P({\bf x})\in S_{\beta
\alpha } $, where $P({\bf x})$ is the projection of point ${\bf
x}$ on $S_{\alpha \beta } $. Let $d_{\,\alpha \beta } ({\bf x})=
\pm \left\| {\bf x}-{\bf x}'\right\| = \pm \sqrt
{\sum_{i=1}^{3}\left( x_i-x'_i \right)^{2}} $ be a distance
between bodies $\Omega _{\alpha } $ and $\Omega _{\beta } $ before
the deformation.

We suppose that possible contact areas $S_{\alpha \beta} $ and $S_{\beta \alpha} $,
$\beta \in B_{\alpha}$, $\alpha = 1,...,N$ \,have nonlinear Winkler covers. Total
compression $w_{\alpha \beta}$ of these covers is related with normal contact stress as
follows: $\sigma _{\alpha  n} ({\bf x})=\sigma _{\beta  n} ({\bf x}')=g_{\alpha
\beta}\left(w_{\alpha \beta}({\bf x})\right)$, ${\bf x} \in S_{\alpha \beta}$, ${\bf
x}' \in S_{\beta \alpha}$, where $g_{\alpha \beta}$ is given nonlinear continuous
function, which satisfy the next conditions:
\begin{equation} \label{prokopyshyni_contrib__6_}
g_{\alpha \beta}(0)=0\,, \,\,\,\,\, \left(\forall \, y,z\right) \,
\left\{\, y<z \, \Rightarrow g_{\alpha \beta}(y) < g_{\alpha
\beta}(z) \,\right\},
\end{equation}
\begin{equation} \label{prokopyshyni_contrib__7_}
\left(\exists \, M_{\alpha \beta}>0 \right) \, \left(\forall \, y,z\right) \, \left\{\,
\left|g_{\alpha \beta}(y)-g_{\alpha \beta}(z)\right| \leq M_{\alpha
\beta}\left|y-z\right| \,\right\}.
\end{equation}

On possible contact zones $S_{\alpha \beta}$, $\beta \in
B_{\alpha}$, $\alpha =1,2,...,N$ we consider the following
unilateral contact conditions through nonlinear Winkler layers:
\begin{equation} \label{prokopyshyni_contrib__8_}
\sigma_{\alpha \xi}({\bf x})=\sigma_{\beta \xi}({\bf x}')=0\,,
\,\,\,\, \sigma_{\alpha \eta}({\bf x})=\sigma_{\beta \eta}({\bf
x}')=0 \,,
\end{equation}
\begin{equation} \label{prokopyshyn_contrib__9_}
\sigma_{\alpha n} ({\bf x})=\sigma_{\beta n} ({\bf x}')= g_{\alpha
\beta}\left( w_{\alpha \beta}({\bf x}) \right)\le 0\, ,
\end{equation}
\begin{equation} \label{prokopyshyn_contrib__10_}
u_{\alpha  n} ({\bf x}) + u_{\beta  n} ({\bf x}') + w_{\alpha
\beta}({\bf x}) \le d_{\,\alpha \beta } ({\bf x}) \, ,
\end{equation}
\begin{equation} \label{prokopyshyni_contrib__11_}
\left[\,u_{\alpha n}({\bf x}) + u_{\beta n}({\bf x}') + w_{\alpha
\beta}({\bf x}) -d_{\alpha \beta } ({\bf x})\,\right]\sigma
_{\alpha n} ({\bf x})=0,\,{\bf x}'=P({\bf x}),\,{\bf x}\in
S_{\alpha \beta }.
\end{equation}

\section{Variational formulations} \label{prokopyshyni_contrib:3}

For each body $\Omega _{\alpha } $ consider Sobolev space
$V_{\alpha } =[H^{1} (\Omega _{\alpha } )]^{3} $ and the closed
subspace $V_{\alpha }^{0} =\left\{{\bf u}_{\,\alpha } \in
V_{\alpha }:\, \, {\bf u}_{\,\alpha } =0 \, \, \, {\rm on} \, \,
\Gamma _{\alpha }^{u} \right\}$. All values of the elements from
these spaces on the parts of boundary $\Gamma _{\alpha }$ should
be understood as traces. The trace of element ${\bf u}_{\alpha}
\in V_{\alpha}$ on the part $\Gamma_{\alpha}^{u}$ should belong to
space $[H^{1/2}(\Gamma_{\alpha}^{u})]^{3}$, and the trace of
element from $V_{\alpha }^{0}$ on the part $\Xi_{\alpha} =
\textrm{int} \,(\Gamma_{\alpha} \setminus \Gamma_{\alpha}^{u})$
should belong to $[H_{00}^{1/2}(\Xi_{\alpha})]^{3}$.

Define Hilbert space $V_{0} = \prod_{\alpha = 1}^{N} V_{\alpha}$
with scalar product \\$\left({\bf u}\,,{\bf v}\right)_{V_{0} }
=\sum _{\alpha =1}^{N}\left({\bf u}_{\,\alpha } ,{\bf v}_{\alpha }
\right)_{V_{\alpha } }  $ and norm $\left\| {\bf u}\right\|
_{V_{0} } =\left({\bf u}\,,{\bf u}\right)_{V_{0}}^{1/2}$, ${\bf
u},{\bf v} \in V_{0}$. Moreover, introduce following spaces $W =
\prod_{\{\alpha,\,\beta\} \in Q} H_{00}^{1/2}({\Xi}_{\alpha}) = \{
{\bf w}=(w_{\alpha \beta})_{\{\alpha,\,\beta\} \in Q}^{\top}:
\,\,\, w_{\alpha \beta} \in H_{00}^{1/2}\}$ and $U_{0} = V_{0}
\times W = \{ {\bf U}=({\bf u},{\bf w})^{\top}:\,\,\,{\bf u} \in
V_{0},\,{\bf w} \in W \}$, where $Q=\left\{\left\{\alpha ,\beta
\right\}:\, \, \alpha \in \left\{1,2,...,N\right\},\, \, \beta \in
B_{\alpha } \right\}$.

In space $U_{0}$ consider the closed convex set of all
displacements, which satisfy nonpenentration contact conditions:
\begin{equation} \label{prokopyshyn_contrib__12_}
K=\{ {\bf U} \in U_{0}:\,\,\, u_{\alpha n} +u_{\beta n}+w_{\alpha
\beta} \le d_{\alpha \beta}\, \, \,{\rm on}\, \, \,S_{\alpha
\beta},\,\,\, \{\alpha,\,\beta\} \in Q \, \},
\end{equation}
where $u_{\alpha \, n}={\bf n}_{\alpha} \cdot {\bf u}_{\alpha} \in
H_{00}^{1/2}(\Xi_{\alpha})$, $w_{\alpha \beta}, d_{\alpha \beta}
\in H_{00}^{1/2}(\Xi_{\alpha})$.

Let us introduce bilinear form $A({\bf u},{\bf v})$, such that
$A({\bf u},{\bf u})$ represents the total elastic deformation
energy of the bodies, linear form $L({\bf u})$, which is equal to
external forces work, and nonquadratic functional $H({\bf w})$,
which represents the total deformation energy of nonlinear Winkler
layers:
\begin{equation} \label{prokopyshyn_contrib__13_}
A\,({\bf u},{\bf v})=\sum _{\alpha=1}^{N}a_{\alpha} ({\bf
u}_{\,\alpha},{\bf v}_{\alpha}),\,\,\,\, a_{\alpha } ({\bf
u}_{\alpha } ,{\bf v}_{\alpha } )= \int _{\Omega _{\, \alpha }
}\hat{\pmb\sigma }_{\alpha } ({\bf u}_{\alpha } )\, \,
:\hat{\pmb\varepsilon }_{\alpha } ({\bf v}_{\alpha } )\,
d\Omega\,,
\end{equation}
\begin{equation} \label{prokopyshyn_contrib__14_}
L\,({\bf u})=\sum _{\alpha=1}^{N}l_{\alpha}({\bf
u}_{\,\alpha}),\,\,\,\, l_{\alpha } ({\bf u}_{\alpha } )= \int
_{\Omega _{\alpha } }{\bf f}_{\,\alpha } \cdot {\bf u}_{\alpha }\,
d\Omega  +\int _{\Gamma _{\alpha }^{\, \sigma } } {\bf p}_{\alpha
} \cdot {\bf u}_{\alpha } \, dS \, \, ,
\end{equation}
\begin{equation} \label{prokopyshyn_contrib__15_}
H({\bf w})=\sum_{\{\alpha,\,\beta\} \in Q} \int_{S_{\alpha \beta}}
\left[ \int_{0}^{w_{\alpha \beta}} g_{\alpha \beta}(z)\,dz \right]
dS\,,\,\,\,\,{\bf u},{\bf v} \in V_{0}\,,\,\,\,\,{\bf w} \in W\,,
\end{equation}
where ${\bf f}_{\,\alpha} \in [L_{2}(\Omega _{\alpha})]^{3}$,
${\bf p}_{\alpha} \in [H_{00}^{-1/2} (\Xi_{\alpha})]^{3}$, $\alpha
= 1,2,...,N$.

We have shown, that if condition (\ref{prokopyshyni_contrib__3_})
holds then bilinear form $A$ is symmetric, continuous and
coercive, and nonquadratic functional $H$ is Gateaux
differentiable:
\begin{equation} \label{prokopyshyni_contrib__16_}
H'({\bf w},{\bf z})=\sum_{\{\alpha,\,\beta\} \in Q}
\int_{S_{\alpha \beta}} g_{\alpha \beta}(w_{\alpha
\beta})\,z_{\alpha \beta}\,dS\,,\,\,\,\, {\bf w}, {\bf z} \in W\,.
\end{equation}

\textbf{Theorem~1.} \textit{Suppose that conditions
(\ref{prokopyshyni_contrib__3_}),
(\ref{prokopyshyni_contrib__6_}), (\ref{prokopyshyni_contrib__7_})
hold. Then problem
(\ref{prokopyshyni_contrib__1_})--(\ref{prokopyshyni_contrib__2_}),
(\ref{prokopyshyni_contrib__4_}),
(\ref{prokopyshyni_contrib__8_})--(\ref{prokopyshyni_contrib__11_})
has an alternative weak formulation as the following minimization
problem: }
\begin{equation} \label{prokopyshyni_contrib__17_}
F({\bf U})=A\,({\bf u},{\bf u})/2-L\,({\bf u})+H({\bf w})\to \mathop{\min
}\limits_{{\bf U}\, \in \, K}.
\end{equation}
\textit{Moreover, there exists a unique solution of problem
(\ref{prokopyshyni_contrib__17_}), and this problem is equivalent
to the following nonlinear variational inequality on set $K$: }
\begin{equation} \label{prokopyshyni_contrib__18_}
F'({\bf U},{\bf V}-{\bf U})=A\,({\bf u},{\bf v}-{\bf u})-L\,({\bf
v}-{\bf u})+H'({\bf w},{\bf z}-{\bf w})\ge 0 \,,\,\,\,\,\forall \,
({\bf v},{\bf z})^{\top} \in K\,.
\end{equation}

Except this variational formulation, we also have proposed another
weak formulation of original contact problem in the form of
nonlinear variational equation.

Let us introduce the following nonquadratic functional in space
$V_{0}$:
\begin{equation} \label {prokopyshyni_contrib__19_}
J({\bf u})=\sum_{\{\alpha,\,\beta\} \in Q} \int_{S_{\alpha \beta}}
\left[ \int_{0}^{d_{\alpha \beta}-u_{\alpha n}- u_{\beta n}}
g_{\alpha \beta}^{-}(z)\,dz \right] dS\,,\,\,\,\,{\bf u} \in
V_{0}\,,
\end{equation}
where $g_{\alpha \beta}^{-}(z)=\{\,0\,,\,\,z \geq 0 \,\} \vee
\{\,g_{\alpha \beta}(z)\,,\,\,z<0\,\}$ is nonlinear function.

Functional $J({\bf u})$ is nonnegative and Gateaux differentiable in $V_{0}$:
\begin{equation} \label {prokopyshyn_contrib__20_}
J'({\bf u},{\bf v})=-\sum_{\{\alpha,\,\beta\} \in Q}
\int_{S_{\alpha \beta}} g_{\alpha \beta}^{-}(d_{\alpha
\beta}-u_{\alpha n}-u_{\beta n})\,\,[v_{\alpha n}+v_{\beta
n}]\,dS\,.
\end{equation}

We have shown that if conditions (\ref{prokopyshyni_contrib__6_})
and (\ref{prokopyshyni_contrib__7_}) hold, then Gateaux
differential $J'({\bf u},{\bf v})$ satisfies the following
properties:
\begin{equation} \label{prokopyshyn_contrib__21_}
(\forall \,{\bf u}\in V_{0})\,\,(\exists \,
\tilde{R}>0\,)\,\,(\forall \, {\bf v}\in V_{0}) \,
\left\{\,\left|J'({\bf u},{\bf v})\right|\le \tilde{R}\left\| {\bf
v}\right\| _{V_{0} } \right\} ,
\end{equation}
\begin{equation} \label{prokopyshyn_contrib__22_}
(\exists \,\tilde{D}>0)\,(\forall \,{\bf u},{\bf v},{\bf w}\!\in\!
V_{0})\! \,\left\{\! \,\left|J'({\bf u}+{\bf w},{\bf
v})\!-\!J'({\bf u},{\bf v})\, \!\right|\!\le\! \tilde{D}\left\|
{\bf v}\right\| _{V_{0}}\! \left\| {\bf w}\right\| _{V_{0}
}\!\,\right\},
\end{equation}
\begin{equation} \label{prokopyshyn_contrib__23_}
\left(\forall \, {\bf u},{\bf v}\in V_{0} \right)\,
\left\{\,J'({\bf u}+{\bf v},{\bf v})-J'({\bf u},{\bf v})\ge 0
\,\right\}.
\end{equation}
These properties helped us to prove the next theorem.

\textbf{Theorem~2.} \textit{Suppose that conditions
(\ref{prokopyshyni_contrib__3_}), (\ref{prokopyshyni_contrib__6_})
and (\ref{prokopyshyni_contrib__7_}) hold. Then the contact
problem
(\ref{prokopyshyni_contrib__1_})--(\ref{prokopyshyni_contrib__2_}),
(\ref{prokopyshyni_contrib__4_}),
(\ref{prokopyshyni_contrib__8_})--(\ref{prokopyshyni_contrib__11_})
is equivalent to problem
(\ref{prokopyshyni_contrib__1_})--(\ref{prokopyshyni_contrib__2_}),
(\ref{prokopyshyni_contrib__4_}), (\ref{prokopyshyni_contrib__8_})
with the following nonlinear boundary value conditions on the
possible contact areas: }
\begin{equation} \label{prokopyshyni_contrib__24_}
\sigma_{\alpha n} ({\bf x})=\sigma_{\beta n} ({\bf x}')= g_{\alpha
\beta}^{-}\left( d_{\alpha \beta}({\bf x})-u_{\alpha n}({\bf
x})-u_{\beta n}({\bf x}') \right)\,,\,\,\,{\bf x}'=P({\bf
x})\,,\,\,\,{\bf x} \in S_{\alpha \beta}\,,
\end{equation}
\textit{and it is equivalent in weak sense to the next
nonquadratic minimization problem: }
\begin{equation} \label{prokopyshyni_contrib__26_}
F_{1}({\bf u})=A\,({\bf u},{\bf u})/2-L\,({\bf u})+J({\bf u})\to
\mathop{\min }\limits_{{\bf u}\, \in \, V_{0} }.
\end{equation}
\textit{Moreover, problem (\ref{prokopyshyni_contrib__26_}) has a
unique solution and is equivalent to the next nonlinear
variational equation in space $V_{0}$: }
\begin{equation} \label{prokopyshyni_contrib__28_}
F'_{1}({\bf u},{\bf v})=A\,({\bf u},{\bf v})+J'({\bf u},{\bf v})-L\,({\bf
v})=0,\,\,\,\, \forall \, {\bf v}\in V_{0}\,, \,\,\, {\bf u}\in V_{0}\,.
\end{equation}

\section{Nonstationary iterative methods}
 \label{prokopyshyni_contrib:5}

In reflexive Banach space $V$ consider an abstract nonlinear
variational equation
\begin{equation} \label{prokopyshyni_contrib__29_}
\Phi \, ({\bf u}, {\bf v})=Y({\bf v})\,, \,\,\,\, \forall \, {\bf v}\in V, \,\,\,\,
{\bf u}\in V,
\end{equation}
where $\Phi:\, \, V \times V \to {\mathbb R}$ is a functional,
which is linear in ${\bf v}$, but nonlinear in ${\bf u}$, and
$Y:\, \, V \to {\mathbb R}$ is linear continuous form. For
numerical solution of (\ref{prokopyshyni_contrib__29_}) consider
the next nonstationary iterative method
\cite{prokopyshyni_contrib_Dyyak2012,prokopyshyni_contrib_Prokopyshyn2012a}:
\begin{equation} \label{prokopyshyni_contrib__30_}
G^{k} ({\bf u}^{k+1} ,{\bf v})=G^{k} ({\bf u}^{k} ,{\bf v})-\gamma ^{k} \left[\,\Phi \,
({\bf u}^{k} ,{\bf v})-Y({\bf v})\,\right], \,\,\,\, k=0,1,... \,\, ,
\end{equation}
where $G^{k}:\, \, V \times V \to {\mathbb R}$ are some given
bilinear forms, ${\gamma}^{k} \in {\mathbb R}$ are iterative
parameters, and ${\bf u}^{k} \in V$ is the \textit{k}-th
approximation to the exact solution of problem
(\ref{prokopyshyni_contrib__29_}).

\textbf{Theorem~3.}\,\cite{prokopyshyni_contrib_Dyyak2012} \textit{\,Suppose that
functional $\Phi$ satisfies the following properties: }
\begin{equation} \label{prokopyshyn_contrib__31_}
\left(\forall {\bf u}\in V \right)\left(\exists R_{\Phi } >0\,
\right)\left(\forall {\bf v}\in V \right)\left\{\, \left|\Phi \,
({\bf u},{\bf v})\right|\le R_{\Phi } \left\| {\bf v}\right\| _{V
} \right\},
\end{equation}
\begin{equation} \label{prokopyshyn_contrib__32_}
\left(\exists D_{\Phi } \!>\!0\right)\!\,\left(\forall {\bf
u},{\bf v},{\bf w}\!\in \!V \right)\!\,\left\{\,\!\,\left|\Phi \,
({\bf u}+{\bf w},{\bf v})-\Phi \, ({\bf u},{\bf v})\,\!\right| \le
D_{\Phi } \!\left\| {\bf v}\right\| _{V } \!\left\| {\bf
w}\right\| _{V } \!\,\right\},
\end{equation}
\begin{equation} \label{prokopyshyn_contrib__33_}
\left(\exists B_{\Phi } >0\right)\left(\forall {\bf u},{\bf v}\in
V \right)\, \left\{\, \Phi \, ({\bf u}+{\bf v},{\bf v})-\Phi \,
({\bf u},{\bf v})\ge B_{\Phi } \left\| {\bf v}\right\| _{V}^{2}
\right\}.
\end{equation}
\textit{Then nonlinear
variational equation (\ref{prokopyshyni_contrib__29_}) has a
unique solution $\bar{{\bf u}} \in V$. In addition, suppose that
bilinear forms} $G^{k}$, $k=0,1,...$ \textit {are symmetric,
continuous with constant} $M_{G}^{*}>0$, \textit {coercive with
constant} $B_{G}^{*}>0$, \textit {and the next conditions hold:}
\begin{equation} \label{prokopyshyn_contrib__34_}
\left(\exists k_{0} \in {\mathbb N}_{0} \right)\left(\forall k\ge
k_{0} \right)\, \left( \forall {\bf u}\in V\right)\, \left\{\,
G^{k} ({\bf u},{\bf u})\ge G^{k+1} ({\bf u},{\bf u})\right\},
\end{equation}
\begin{equation} \label{prokopyshyn_contrib__35_}
\left(\exists \varepsilon \in (0,\, \gamma ^{*} ),\,\, \gamma^{*}
={B_{\Phi } B_{G}^{*} \mathord{\left/ {\vphantom {B_{\Phi }
B_{G}^{*}  D_{\Phi }^{2} }} \right. \kern-\nulldelimiterspace}
D_{\Phi }^{2} } \right)\left(\exists k_{1} \right)\left(\forall
k\ge k_{1} \right)\, \left\{\, \gamma ^{k} \in [\varepsilon , \,
2\gamma ^{*} -\varepsilon ]\, \right\},
\end{equation}
\textit{where $\{ {\bf u}^{k} \} \subset V$ is obtained by
iterative method (\ref{prokopyshyni_contrib__30_}). }

\section{Domain decomposition schemes} \label{prokopyshyni_contrib:5}

Now let us apply nonstationary iterative method
(\ref{prokopyshyni_contrib__30_}) for solving nonlinear
variational equation (\ref{prokopyshyni_contrib__28_}), which
corresponds to original contact problem. This equation can be
written in form (\ref{prokopyshyni_contrib__29_}), where
$\Phi({\bf u},{\bf v}) = A\,({\bf u},{\bf v})+J'({\bf u},{\bf
v})\,$, $Y({\bf v})=L\,({\bf v})\,$, ${\bf u},{\bf v} \in V$,
$V=V_{0}\,$, and iterative method
(\ref{prokopyshyni_contrib__30_}) applied to solve
(\ref{prokopyshyni_contrib__28_}) rewrites as follows:
\begin{equation} \label{prokopyshyni_contrib__37_}
G^{k}({\bf u}^{k+1}, {\bf v})=G^{k}({\bf u}^{k}, {\bf v})-\gamma
^{k} \left[A\,({\bf u}^{k}, {\bf v})+J'({\bf u}^{k}, {\bf
v})-L({\bf v}) \right],\,k=0,1,....
\end{equation}

Note, that in general case iterative method
(\ref{prokopyshyni_contrib__37_}) does not lead to domain
decomposition. Let us propose such variants of this method, which
involve the domain decomposition. At first, let us take bilinear
forms $G^{k}$ in method (\ref{prokopyshyni_contrib__37_}) as
follows:
\begin{equation} \label{prokopyshyni_contrib__38_}
G^{k}({\bf u},{\bf v})=\partial^{2} F_{1}({\bf u}^{k},{\bf u},{\bf v})=A\,({\bf u},{\bf
v})+\partial^{2} J({\bf u}^{k},{\bf u},{\bf v})\,,\,\,\,\, {\bf u},{\bf v} \in V_{0}\,,
\end{equation}
$$
\partial^{2} J({\bf u}^{k},{\bf u},{\bf v})=\sum_{\{\alpha,\,\beta \}\in Q}
\int_{S_{\alpha \beta } }\chi _{\alpha \beta }^{k} \, g'_{\alpha \beta } (d_{ \alpha
\beta } -u_{\alpha n}^{k} -u_{\beta n}^{k} )\left[u_{\alpha n} +u_{\beta n}
\right]\left[v_{\alpha n} +v_{\beta n} \right]dS,
$$
\begin{equation} \label{prokopyshyni_contrib__39_}
{\chi}_{\alpha \beta}^{k} = -[\,\textrm{sgn}\,(d_{\alpha \beta }
-u_{\alpha n}^{k} -u_{\beta n}^{k})\,]^{-}=\{\,0,\, d_{\alpha
\beta } -u_{\alpha n}^{k} -u_{\beta n}^{k} \ge 0 \,\}\vee \{\,1,\,
\textrm{else}\,\}.
\end{equation}
Here $\partial^{2} F_{1}({\bf u}^{k},{\bf u},{\bf v})$ and
$\partial^{2} J({\bf u}^{k},{\bf u},{\bf v})$ are the second
subdifferentials of functionals $F_{1}$ and $J$ in point ${\bf
u}^{k} \in V_{0}$. In the case when ${\gamma}^{k}=1$,
$k=0,1,...\,$, iterative method (\ref{prokopyshyni_contrib__37_})
with bilinear forms (\ref{prokopyshyni_contrib__38_}) corresponds
to semismooth Newton method for variational equation
(\ref{prokopyshyni_contrib__28_}). However, this method does not
lead to domain decomposition.

Now, let us take bilinear forms $G^{k}$ in the following way:
\begin{equation} \label{prokopyshyni_contrib__40_}
G^{k}({\bf u},{\bf v})=A\,({\bf u},{\bf v})+X^{k}({\bf u},{\bf v})\,,\,\,\,\, {\bf
u},{\bf v}\in V_{0}\,,
\end{equation}
\begin{equation} \label{prokopyshyni_contrib__40b_}
X^{k}({\bf u},{\bf v})=\sum_{\alpha =1}^{N}\sum _{\,\beta \in
B_{\alpha } }\int_{S_{\alpha \beta}} \psi_{\alpha \beta}^{k} \,
g'_{\alpha \beta } (d_{\alpha \beta } -u_{\alpha n}^{k} -u_{\beta
n}^{k} )\, u_{\alpha n} v_{\alpha n}  dS, \,\,{\bf u},{\bf v} \in
V_{0},
\end{equation}
where $\psi_{\alpha \beta}^{k}({\bf x}) = \{\,1,\, {\bf x}\in
S_{\alpha \beta}^{k} \,\}\vee \{\,0,\, {\bf x}\in S_{\alpha \beta
} \backslash S_{\alpha \beta }^{k} \,\}$ are characteristic
functions of some given subsets $S_{\alpha \beta}^{k} \subseteq
S_{\alpha \beta}$ of possible contact areas.

Let us show, that such choice of bilinear forms $G^{k}$ involves
the domain decomposition. Introduce a notation $\tilde{{\bf
u}}^{k+1}=({\bf u}^{k+1} -{\bf u}^{k})/{\gamma}^{k}+{\bf u}^{k}
\in V_{0}$. Then iterative method
(\ref{prokopyshyni_contrib__37_}) with bilinear forms
(\ref{prokopyshyni_contrib__40_}) can be written in such way:
\begin{equation} \label{prokopyshyni_contrib__41_}
A\,(\tilde{{\bf u}}^{k+1},{\bf v})+X^{k}(\tilde{{\bf
u}}^{k+1},{\bf v})=L\,({\bf v})+X^{k}({\bf u}^{k},{\bf v})-J'({\bf
u}^{k},{\bf v})\,,\,\,\,\forall\,{\bf v} \in V_{0}\,.
\end{equation}
\begin{equation} \label{prokopyshyni_contrib__42_}
{\bf u}^{k+1} ={\gamma}^{k}\,\tilde{{\bf u}}^{k+1} +
(1-{\gamma}^{k})\,{\bf u}^{k}, \,\,\,\, k=0,1,...\,.
\end{equation}

Since the common quantities of the subdomains are known from the
previous iteration, variational equation
(\ref{prokopyshyni_contrib__41_}) splits into $N$ separate
equations in subdomains $\Omega_{\alpha}\,$, and iterative method
(\ref{prokopyshyni_contrib__41_})--(\ref{prokopyshyni_contrib__42_})
can be written in the following equivalent form:
\[a_{\alpha}(\tilde{{\bf u}}_{\alpha }^{k+1},
{\bf v}_{\alpha }) + \sum_{\beta \, \in B_{\alpha } }\int _{S_{\alpha \beta } }\psi
_{\alpha \beta }^{k} \, g'_{\alpha \beta } (d_{ \alpha \beta } -u_{\alpha n}^{k}
-u_{\beta n}^{k} )\, \tilde{u}_{\alpha n}^{k+1} v_{\alpha n} \, dS = \]
\[=l_{\alpha}({\bf v}_{\alpha})+\sum _{\beta \, \in B_{\alpha } }\int _{S_{\alpha
\beta } }\psi _{\alpha \beta }^{k} \, g'_{\alpha \beta } (d_{\alpha \beta } -u_{\alpha
n}^{k} -u_{\beta n}^{k} )\, u_{ \alpha n}^{k} v_{\alpha n}\,dS \, + \]
\begin{equation} \label{prokopyshyni_contrib__43_}
+\sum _{\beta \, \in B_{\alpha } }\int _{S_{\alpha \beta } }g_{\alpha \beta }^{-}
(d_{\alpha \beta } -u_{\alpha n}^{k} -u_{\beta n}^{k} )\, v_{\alpha n}\,dS\,,\,\,\,\,
\forall\,{\bf v}_{\alpha} \in V_{\alpha}^{0}\,,
\end{equation}
\begin{equation} \label{prokopyshyni_contrib__44_}
{\bf u}_{\alpha}^{k+1} = {\gamma}^{k}\,\tilde{{\bf
u}}_{\alpha}^{k+1} + (1-{\gamma}^{k})\,{\bf u}_{\alpha }^{k}\,,
\,\,\,\, \alpha =1,2,...,N, \,\,\,\, k=0,1,...\,.
\end{equation}

In each iteration $k$ of method
(\ref{prokopyshyni_contrib__43_})--(\ref{prokopyshyni_contrib__44_}),
we have to solve $N$ linear variational equations
(\ref{prokopyshyni_contrib__43_}) in parallel, which correspond to
linear elasticity problems in separate bodies $\Omega_{\alpha}$
with Robin boundary conditions on possible contact areas.
Therefore, this method refers to parallel Robin--Robin type domain
decomposition schemes.

By taking different characteristic functions $\psi_{\alpha
\beta}^{k}$, we can obtain different particular cases of domain
decomposition method
(\ref{prokopyshyni_contrib__43_})--(\ref{prokopyshyni_contrib__44_}).
Thus, taking $\psi_{\alpha \beta}^{k}({\bf x}) \equiv 0$
$(S_{\alpha \beta}^{k} = \emptyset)$, $\forall \alpha, \beta$,
$\forall k$, we get parallel Neumann--Neumann domain decomposition
scheme. Other borderline case is when $\psi_{\alpha
\beta}^{k}({\bf x}) \equiv 1$ $(S_{\alpha \beta}^{k} = S_{\alpha
\beta})$, $\forall \alpha, \beta$, $\forall k$.

Moreover, we can choose characteristic functions $\psi_{\alpha
\beta}^{k}$ by formula (\ref{prokopyshyni_contrib__39_}), i.e.
$\psi_{\alpha \beta}^{k}=\chi_{\alpha \beta}^{k}$. Numerical
experiments, provided by us, have shown, that such DDM has higher
convergence rate than other particular domain decomposition
schemes.

\section{Numerical investigations} \label{prokopyshyni_contrib:6}

Numerical investigations of proposed DDMs have been made for plane
problem of unilateral contact between two isotropic bodies
${\Omega}_{1}$ and ${\Omega}_{2}$, one of which has a groove
(Fig.~1b). The bodies are uniformly loaded by normal stress with
intencity $q=10\,\textrm{MPa}$. Each body has length
$l=4\,\textrm{cm}$ and height $h=1\,\textrm{cm}$. The Young's
moduli and Poisson's ratios of the bodies are the same:
$E_{1}=E_{2}=2.1 \cdot 10^{5}\,\textrm{MPa}$,
${\nu}_{1}={\nu}_{2}=0.3$. The distance between bodies is $d_{12}
({\bf x})=r\left\{{\, [\, 1-(x_{1} -l)^{2} \mathord{\left/
{\vphantom {\, [\, 1-(x_{1} -l)^{2} b^{2} }} \right.
\kern-\nulldelimiterspace} b^{2} } ]^{+} \right\}^{3/2}$, ${\bf x}
\in S_{12}\,$, where $b=1\,\textrm{cm}$, $r=5 \cdot
10^{-4}\,\textrm{cm}$, $z^{+}=\textrm{max}\,\{0,z\}$, $S_{12}
=\left\{{\bf x}=(x_{1},x_{2})^{\top}:\,\,\, x_{1} \in [0,\, l],\,
\, x_{2} =h\right\}$.

Across possible contact area $S_{12}$ there is a nonlinear Winkler
layer. The relationship between normal contact stresses and
displacements of this layer are described by the following power
function:\\ $g_{12}\left(w_{12}({\bf
x})\right)=B^{-1/a}\,\textrm{sgn}\left(w_{12}({\bf x})\right)
\left|w_{12}({\bf x})\right|^{1/a}$, ${\bf x} \in S_{12}\,$, where
parameters $B$ and $a$ are taken from the intervals $B \in
[\,10^{-6}\,\textrm{cm}/(\textrm{MPa})^{a},$ $2 \cdot
10^{-4}\,\textrm{cm}/(\textrm{MPa})^{a}\,]\,$, $a \in
[\,0.1,\,1\,]$. For such choice of these parameters the nonlinear
Winkler layer models a roughness of the possible contact surface
\cite{prokopyshyni_contrib_Goryacheva1998}.

This problem has been solved by DDM
(\ref{prokopyshyni_contrib__43_})--(\ref{prokopyshyni_contrib__44_})
with stationary iterative parameters ${\gamma}^{k}=\gamma$,
$\forall \, k$ and characteristic functions $\psi_{12}^{k}$, taken
by formula (\ref{prokopyshyni_contrib__39_}), i.e.
$\psi_{12}^{k}=\chi_{12}^{k}$, $\forall \, k$. For solving linear
variational problems (\ref{prokopyshyni_contrib__43_}) in each
iteration $k$ we have used finite element method with 8192 linear
triangular elements for each body.

We have used the following initial guesses for displacements $u_{1
n}^{0}({\bf x}) \equiv {10}^{-4}$, $u_{2 n}^{0}({\bf x}) \equiv
{10}^{-4}$ and the next stopping criterion:
${\rho}_{\alpha}^{k+1}=\left\| u_{\alpha n}^{k+1}-u_{\alpha n}^{k}
\right\|_{2}\,/$ $\left\| u_{\alpha n}^{k+1} \right\|_{2} \le
{\varepsilon}_{u}$, $\alpha = 1,2$, where $\left\| u_{\alpha \,n}
\right\|_{2} =\sqrt{\sum _{j}\left[u_{\alpha \,n}({\bf
x}^{j})\right]^{2}}$ is discrete norm, ${\bf x}^{j} \in S_{12}$
are finite element nodes on the possible contact area, and
$\varepsilon_{u}>0$ is relative accuracy.

\begin{figure}[h]
%\sidecaption
\includegraphics[scale=0.51]{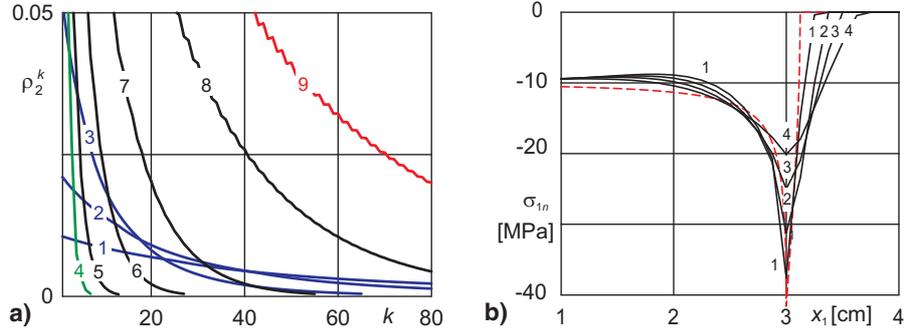}
\caption{Relative error (a), and normal contact stress (b) }
\label{author_program_element_fig:2}       % Give a unique label
\end{figure}
At Fig.~2a the relative error ${\rho}_{2}^{k}$ of displacement
$u_{2 n}$ on different iterations $k$, obtained for $B=2.5 \cdot
{10}^{-5}\,\textrm{cm}/(\textrm{MPa})^{a}$, $a=0.5$, is
represented for different values of parameter $\gamma$.
Curves~1--9 correspond to $\gamma=0.01$, 0.02, 0.05, 0.6, 0.8
(0.3), 0.9, 0.95, 0.98, 0.99. For these values of parameter
$\gamma$, DDM
(\ref{prokopyshyni_contrib__43_})--(\ref{prokopyshyni_contrib__44_})
reaches the accuracy $\varepsilon_{u}=10^{-3}$ in 193, 124, 65, 7,
13, 27, 55, 134 iterations respectively.

Thus, we conclude, that the best convergence rate reaches if
$\gamma=0.6$. The convergence rate is good if $\gamma \in
[\,0.1,\,0.9\,]$. However, it becomes slow when $\gamma$ is close
to 0 or to 1. For $\gamma=0.98$ the method is still convergent,
but the convergence becomes nonmonotone. For $\gamma \geq 0.99$
the method is not anymore convergent.

At Fig.~2b the normal contact stress ${\sigma}_{1 n}={\sigma}_{2
n}$, obtained by DDM
(\ref{prokopyshyni_contrib__43_})--(\ref{prokopyshyni_contrib__44_})
for $B=10^{-5}\,\textrm{cm}/(\textrm{MPa})^{a}$ and different
values of parameters $a$ is represented. Curves~1--4 correspond to
numerical solution for $a=0.3$, 0.6, 0.8, 1. Dashed curve
represents the analytical solution, obtained in
\cite{prokopyshyni_contrib_Shvets1996} for %the case of unilateral
contact of two halfspaces without nonlinear layer. Here we
conclude, that for small values of $a$ ($a \leq 0.3$) the
influence of nonlinear %Winkler
layer on the contact behavior is not so large and the numerical
solutions are close to the solution without layer. However, for
larger values of $a$ ($a \geq 0.5$) the influence of nonlinear
layer becomes more significant and can not be neglected.

The positive feature of proposed domain decomposition methods are
the simplicity of their algorithms. These methods have only one
iteration loop, which deals with domain decomposition,
nonlinearity of Winkler layers and contact conditions.

%

%% BIBLIOGRAPHY %%%%%%%%%%%%%%%%%%%%%%%%%%%%%%%%%%%%%%%%%%%%%

%\bibliography{prokopyshyni_contrib}

\end{document}